\documentclass[12pt]{elsarticle}
\usepackage{txfonts}
\usepackage{graphicx}
\usepackage{tikz}

\usetikzlibrary{positioning}

\newcommand{\n}{\noindent}
\newtheorem{thm}{Theorem}

\newtheorem{lem}{Lemma}

\newtheorem{transformation}{Transformation}
\newtheorem{remark}{Remark}

\begin{document}

\title{Ordering Unicyclic Graphs with Respect to F-index}

\author[man]{Ruhul Amin}
\ead{aminruhul80@gmail.com}
\author[man]{Sk. Md. Abu Nayeem\corref{cor1}}
\ead{nayeem.math@aliah.ac.in}

\address[man]{Department of Mathematics, Aliah University, New Town, Kolkata -- 700 156, India.}
\cortext[cor1]{Corresponding Author.}
\begin{abstract}
F-index of a graph is the sum of the cube of the degrees of the vertices. In this paper, we investigate the F-indices of unicyclic graphs by introducing some transformation, and characterize the unicyclic graphs with the first five largest F-indices and the unicyclic graphs with the first two smallest F-indices, respectively.
\end{abstract}

\maketitle
\section{Introduction}
A topological index is a numerical value associated with graph representation of a molecule for correlation of chemical structure of the molecule with various physical properties, chemical reactivity or biological activity. The first Zagreb index, introduced in 1972 \cite{gut04}, is one of the oldest topological indices. It is defined as $$ M_{1}(G)= \displaystyle\sum_{v \in V(G)}d^2(v) = \displaystyle\sum_{uv \in E(G)}[d(u) + d(v)].$$

In the same paper where first Zagreb index was introduced by Gutman and Trinajsti\'{c} \cite{gut72}, it was shown that the sum of squares and sum of the cubes of the vertex degrees of the underlying molecular graph influence the total $\pi$-electron energy $E$. The sum of squares of the vertex degrees being known as the first Zagreb index, has been subject of hundreds of researches, but the sum of cubes of the vertex degrees remained unstudied for a long time by scholars doing research on degree--based topological indices. Furtula and Gutman, restudied this index recently and named it as ``forgotten'' topological index, or $F$-index \cite{fur15}. Thus, $F$-index is defined as $$F(G)=\sum\limits_{v\in V(G)}{{{d}_{G}}{{(v)}^{3}}}=\sum\limits_{uv\in E(G)}{[{{d}_{G}}{{(u)}^{2}}+{{d}_{G}}{{(v)}^{2}}]}.$$

A leaf is a vertex of degree one and pendant edges are edges incident to a leaf and a stem, where a stem is a vertex adjacent to at least one leaf. $P_{n}, C_{n}$ and $K_{1, n - 1}$ respectively denotes path, cycle, and star with $n$ vertices. $\mathcal{U}_{n}$ and $\mathcal{U}_{n}^{k}$ denotes the set of all unicyclic graph with $n$ vertices and the set of all unicyclic graph with $n$ vertices and cycle length $k$ respectively. We denote  the unicyclic graph constructed by attaching $n - k$ leaves to one vertex on a cycle of length $k$ by $G_{K, 1}^{(n)}$ as shown in Figure 1(a). $G_{K, 2}^{(n)}$ denotes the unicyclic graph constructed by attaching $n - k - 1$ leaves to one vertex $u$ of the cycle and $K_{2}$ attached to adjacent vertex of $u$ as shown in Figure 1(b). We denote the unicyclic graph  constructed by attaching $K_{1, p}, K_{1, q}$ and $K_{1, r}$ to the vertices of $C_{3}$ respectively by $S_{p, q, r}$ $(p, q, r \geq 0, p + q + r =n - 3)$ as shown in Figure 1(c).

Extremal Zagreb indices of unicyclic graphs has been studied by Zhou \cite{zho}. Ordering of unicyclic graphs with respect to Zagreb indices has been studied by Xia \cite{fan07}. Unicyclic graphs with the  first three smallest and largest first general Zagreb index has been studied by Zhang \cite{zha06}. A unified approach to the extremal Zagreb indices for trees, unicyclic graphs, and  bicyclic graphs has been proposed by Deng \cite{den07}. Extremal trees with respect to F-index have been found by Abdo et al.\cite{abd15}. F-index of some graph operations has been studied by De et al.\cite{de16}.

In this paper, we investigate the F-indices of unicyclic graphs by introducing some transformation, and characterize the unicyclic graphs for the first five largest F-indices and the unicyclic graphs with the  first two smallest F-indices.

\section{Two transformations which increase the F-indices}
Let $ E_{1} \subseteq E(G)$. We denote by $G - E_{1}$ the subgraph of $G$ obtained by deleting the edges in $E_{1}$. Let $W \subseteq V(G)$. $G - W$ denotes the subgraph of $G$ obtained by deleting the vertices in $W$ and the edges incident with them. Again let, $E_2 \subseteq E(\overline{G})$, where $\overline{G}$ is the complement of $G$. Then by $G + E_2$ we mean the graph obtained by adding the edges in $E_2$ to $G$.

We give two transformations which will increase the F-indices as follows.
\begin{transformation}\label{A}
Let $uv$ be an edge of the graph $G$, $d_{G}(u) \geq 2$, $N_{G}(v) = \{u, w_{1}, w_{2}, \ldots, w_{s}\}$, and $w_{1}, w_{2}, \ldots, w_{s}$ are leaves. Then $G_{1} = G - \{vw_{1}, vw_{2}, \ldots, vw_{s}\} + \{uw_{1}, uw_{2}, \ldots, uw_{s}\}$, as shown in Figure 2, is said to be the graph obtained from $G$ by Transformation A.
\end{transformation}

\begin{lem}\label{21}
Let $G_{1}$ be obtained from $G$ by Transformation $A$, then $F(G_{1}) > F(G)$.
\end{lem}

\noindent\textit{Proof.}
Let $G_{0} = G - \{v, w_{1}, w_{2}, \ldots, w_{s}\}$.

By definition of F-index, we have
\begin{eqnarray*}
F(G_{1}) - F(G)& =& d^3_{G_{1}}(v) - d^3_{G}(v) + d^3_{G_{1}}(u) - d^3_{G}(u)\\
&=&( d^3_{G}(v) + s)^3 - d^3_{G}(v) + 1 - (s + 1)^3\\
&=&3(sd^2_{G}(v) - s^2) + 3s (d_{G}(v) - 1)\\
&>&0 ~(\mbox{since}~ d_G(v) > s).
\end{eqnarray*}

Hence, $F(G_{1}) > F(G)$.
\qed

\begin{remark}\label{1}
Using Transformation A repeatedly, any unicyclic graph can be transformed into such a unicyclic graph that all the edges not on the cycle are pendant edges.
\end{remark}

\begin{transformation}\label{B}
Let $u$ and $v$  be two vertices in $G$. $u_{1}, u_{2}, \ldots, u_{s}$ are the leaves adjacent to $u$, $v_{1}, v_{2}, \ldots, v_{t}$ are the leaves adjacent to $v$. Then $G_{1} = G - \{uu_{1}, uu_{2}, \ldots, uu_{s}\} + \{vu_{1}, vu_{2}, \ldots, vu_{s}\}$, and $G_{2} = G - \{vu_{1}, vu_{2}, \ldots, vu_{t}\} + \{uv_{1}, uv_{2}, \ldots, uv_{t}\}$, as shown in Figure 3, are said to be the graphs obtained from $G$ by Transformation B.
\end{transformation}

\begin{lem}\label{22}
Let $G_{1}$ and $G_{2}$ be obtained from $G$ by Transformation $B$. Then either $F(G_{1}) > F(G)$ or $F(G_{2}) > F(G)$.
\end{lem}

\noindent\textit{Proof.}
Let $G_{0} = G - \{u_{1}, u_{2}, \ldots, u_{s}, v_{1}, v_{2}, \ldots, v_{t}\}$.

Then by the definition of F-index, we have
\begin{eqnarray*}
F(G_{1}) - F(G)& = &d^3_{G_{1}}(v) - d^3_{G}(v) + d^3_{G_{2}}(u) - d^3_{G}(u)\\
&=& (d^3_{G}(v) + s) - d^3_{G}(v) + (d^3_{G}(u) - s) - d^3_{G}(u)\\
&=& 3s(d_{G}(v) + d_{G}(u)) (d_{G}(v) - d_{G}(u) + s).\\
F(G_{2}) - F(G)& = &d^3_{G_{2}}(v) - d^3_{G}(v) + d^3_{G_{2}}(u) - d^3_{G}(u)\\
& =& (d^3_{G_{2}}(v) -t)^3 - d^3_{G}(v) + (d^3_{G_{2}}(u) + t)^3 - d^3_{G}(u)\\
&=& 3t (d_{G}(u) + d_{G}(v)) (d_{G}(u) - d_{G}(v) + t).
\end{eqnarray*}
Since either $(d_{G}(v) - d_{G}(u) + s) > 0$ or $ (d_{G}(u) - d_{G}(v) + t) > 0$,
it is clear that either $F(G_{1}) > F(G)$ or $F(G_{2}) > F(G)$.
\qed

\begin{remark}\label{2}
Using Transformation B repeatedly, any unicyclic graph can be transformed into such a unicyclic graph that all the pendant edges are attached to the same vertex.
\end{remark}

\section{Some transformations which decrease the F-indices}
\begin{transformation}\label{C}
Let $G \neq P_{1}$ be a connected graph and we choose $u\in V(G)$ and $G_{1}$ denotes the graph that results from identifying $u$ with the vertex $v_{k}$ of a simple path $v_{1}v_{2}...v_{n}, 1 < k < n$. The graph $G_{2}$  obtained from $G_{1}$ by deleting $v_{k-1}v_{k}$ and adding $v_{k-1}v_{n}$, as shown in Figure 4, is said to be the graph obtained by applying Transformation C to the graph $G_1$.
\end{transformation}

\begin{lem}\label{31}
Let $G_{2}$ be the graph obtained by Transformation C to $G_{1}$, as shown in Figure 4. Then $F(G_{1}) > F(G_{2})$.
\end{lem}

\noindent\textit{Proof.}
By the definition of F-index, we have
\begin{eqnarray*}
F(G_{1}) - F(G_{2}) &=& 1^3 + 1^3 - 1^3 + d_{G_{1}}^3(v_{k}) - d_{G_{2}}^3(v_{k}) - 2^3\\
&=& 1 + d_{G_{1}}^3(v_{k}) - (d_{G_{1}}(v_{k}) - 1)^3 - 2^3\\
&=&3d_{G_{1}}^2(v_{k}) - 3d_{G_{1}}(v_{k}) -  6\\
 &>& 0~[\mbox{since}~d_{G_1}(v_k) = 4].
\end{eqnarray*}
Hence, $F(G_{1}) > F(G_{2})$.
\qed

\begin{remark}\label{3}
By Repeated application of Transformation C, any tree $T$ attached to a graph $G$ can be transformed into a path attached to $G$ and F-index will decreases at each step of applying Transformation C.
\end{remark}

\begin{transformation}\label{D}
Let $u$ and $v$  be two vertices in a graph $G$. $G_{1}$ denotes the graph that results from identifying $u$ with the vertex $u_{0}$ of a path $u_{0}u_{1}u_{2}...u_{s}$ and identifying $v$ with the vertex $v_{0}$ of a path $v_{0}v_{1}v_{2}...v_{t}$. The graph $G_{2}$ obtained from $G_{1}$ by deleting $uu_{1}$ and adding $v_{t}u_{1}$, as shown in Figure 5, is the graph obtained by applying Transformation D to $G_1$.
\end{transformation}

\begin{lem}\label{32}
Let $G_{2}$ be the graph obtained from $G_1$, by Transformation D, as shown in Figure 5, where $d_{G}(u) \geq d_{G}(v) > 1, s \geq 1$ and $t\geq 0. $\\
(i) If $t > 0$, then $F(G_{1}) > F(G_{2})$.\\
(ii) If $t = 0$, and $d_G(u) > d_G(v)$, then $F(G_{1}) > F(G_{2})$.\\
\end{lem}

\n\textit{Proof.} (i) Since $d_{G}(u) > 1$ and $t > 0$, we have
\begin{eqnarray*}
F(G_{1}) - F(G_{2})& =& d^3_{G_{1}}(u) + d^3_{G_{1}}(v_{t}) - d^3_{G_{2}}(u) - d^3_{G_{2}}(v_{t})\\
&=& (d_{G}(u) + 1)^3 + 1 - d^3_{G}(u) - 8\\
&=& 3d^2_{G}(u) + 3d_{G}(u) - 6 > 0.
\end{eqnarray*}
Hence, $F(G_{1}) > F(G_{2})$.

\noindent (ii) If $t = 0$ and $d_{G}(u) > d_{G}(v)$ then
\begin{eqnarray*}
F(G_{1}) - F(G_{2}) &=& d^3_{G_{1}}(u) + d^3_{G_{1}}(v_{t}) - d^3_{G_{2}}(u) - d^3_{G_{2}}(v_{t})\\
&=& (d_{G}(u) + 1)^3 + d^3_{G}(v) - d^3_{G}(u) - d^3_{G}(v)\\
&=& 3d^2_{G}(u) + 3d_{G}(u) + 1 > 0.
\end{eqnarray*}
Hence, $F(G_{1}) > F(G_{2})$.
\qed

\begin{remark}
After repeated application of Transformation D, any tree attached to a unicyclic graph can be transformed into such a unicyclic graph that a path is attached to a cycle, and F-indices decrease in each application of Transformed D.
\end{remark}

\begin{transformation}\label{E}
Let $G_{1}$ be a unicyclic graph constructed by attaching $n - k$ leaves to a vertex $u$ on a cycle of length $k$. The graph $G_{2}$ obtained from $G_{1}$ by attaching a path of length $n - k$ at $u$, as shown in Figure 6, is said to be the graph obtained from $G_1$ by Transformation E.
\end{transformation}
\begin{lem}\label{33}
Let the graph $G_{2}$ is obtained from $G_{1}$ by Transformation $E$. Then $F(G_{2}) \leq F(G_{1})$.
\end{lem}

\n\textit{Proof.}
By the definition of F-index, we have\\
$F(G_{1}) - F(G_{2}) = (n - k + 2)^3 + (n - k).1^3 - 3^3 - (n - k 1).2^3 - 1$ $$= (n - k)^3 + 6(n - k)^2 + 5(n - k) - 12 \geq 0 ~\mbox{since}~(n - k)\geq 1.$$

Hence, $F(G_{2}) \leq F(G_{1})$.
\qed

\begin{transformation}\label{F}
Let $G_{1}$ be a unicyclic graph constructed by attaching $n - k$ leaves to a vertex $u$ on a cycle of length $k$. By Transformation F, the graph $G_{2}$ is obtained from $G_{1}$ by attaching two paths of length $s$ and $t$, where $s + t = n - k$, at the vertex $u$, as shown in Figure 7.
\end{transformation}
\begin{lem}\label{34}
Let the graph $G_{2}$ is obtained from $G_{1}$ by Transformation $F$. Then $F(G_{2}) \leq F(G_{1})$.
\end{lem}

\n\textit{Proof.}
By the definition of F-index, we have\\
$F(G_{1}) - F(G_{2}) = (n - k + 2)^3 + (n - k).1^3 - 4^3 - (s - 1)2^3 - (t - 1)2^3 - 2.1^3$ $$= (n - k)^3 + 6(n - k)^2 + 5(n - k) - 42 \geq 0,~ \mbox{since} ~(n - k) \geq 3.$$
Hence, $F(G_{2}) \leq F(G_{1})$.
\qed

\section{Unicyclic graphs with larger F-indices}
In this section, we obtain some upper bounds of the unicyclic graphs with respect to their F-indices.
\begin{lem}\label{41}
 Let $G \in \mathcal{U}_n^{k}$. Then $F(G) \leq F(G_{k ,1}^{(n)})$, with equality if and only if $G \cong G_{k ,1}^{(n)}$.
\end{lem}

\n\textit{Proof.}
Using Lemma 1 to Lemma 6, we have $F(G) \leq F(G_{k, 1}^{(n)})$ and equality holds if and only if $G \cong G_{k, 1}^{(n)}$.
\qed

\begin{lem}\label{42}
 Let $G \in \mathcal{U}_n$. Then\\
 (i) $F(G) \leq F(G_{3, 1}^{(n)})$, with equality if and only if $G \cong G_{3, 1}^{(n)}$.\\
 (ii) If $G$ is not isomorphic to $G_{3, 1}^{(n)}$, then $F(G) \leq F(G_{3, 2}^{(n)})$ with equality if and only if $G \cong G_{3, 2}^{(n)}$.
 \end{lem}
 
 \n\textit{Proof.} (i)
 By Lemma 7, we have $F(G) \leq F(G_{k ,1}^{(n)})$.

  Next we prove that $F(G_{k, 1}^{(n)}) \leq F(G_{3, 1}^{(n)})$.

   By definition of F-index, we have $$ F(G_{3, 1}^{(n)}) - F(G_{k ,1}^{(n)}) = (n - 1)^3 - (n - k + 2)^3 + 2.2^3 - 2^3.(k - 1) + (n - 3) - (n - k)$$$$ = (k - 3)\{(n - 1)^2 + (n - 1)(n - k + 2) + (n -k + 2)^2 - 7\} \geq 0 \mbox{~since}~ n - k \geq 1 \mbox{ and} ~k \geq 3.$$

 Hence, $F(G_{k, 1}^{(n)}) \leq F(G_{3, 1}^{(n)})$.  Therefore, $F(G) \leq F(G_{3, 1}^{(n)})$, with equality if and only if $G \cong G_{3, 1}^{(n)}$.

\noindent(ii)
 By Lemma 7, we have, $F(G) \leq F(G_{k ,1}^{(n)})$.

  Next we prove that if $G$ is not isomorphic to $G_{3, 1}^{(n)}$, then $ G_{k, 1}^{(n)} \leq  G_{3, 2}^{(n)}$.

   By definition of F-index, we have $$F(G_{3, 2}^{(n)}) - F(G_{k, 1}^{(n)}) = (n - 2)^3 + 3^3 + 2^3 + (n - 3) - (n - k + 2)^3 - (k - 1).2^3 - (n -k)$$ $$= (k - 4)\{3n(n - k) + (k -2)^2\} - 8(k - 4) + 8 > 0~ if ~k \geq 4.$$

If $k = 3$ and $G_{k, 1}^{(n)}$ is not isomorphic to $G_{3, 1}^{(n)}$, then obviously $G_{k, 1}^{(n)} \leq G_{3, 2}^{(n)}$.

Hence, $F(G) \leq F(G_{3, 2}^{(n)})$, with equality if and only if $G \cong G_{3, 2}^{(n)}$.
\qed

\begin{thm}\label{41}
 Let $G \in \mathcal{U}_n^{k}(k \geq 3)$. Then $F(G_{k, 1}^{(n)}) > F(G_{k, 2}^{(n)})$.
\end{thm}

\n\textit{Proof.}
From the definition of F-index, we have
\begin{eqnarray*}
F(G_{k ,1}^{(n)}) - F(G_{k ,2}^{(n)}) &=& (n - k - 2)^3 + 2^3 - 3^3 - (n - k + 1)^3\\
&=& 3(n - k)^2 + 9(n - k) - 2 > 0 ~\mbox{since} ~n - k > 0.
\end{eqnarray*}

Hence, $F(G_{k ,1}^{(n)}) > F(G_{k ,2}^{(n)})$.
\qed

\begin{thm}\label{42}
 Let $G \in \mathcal{U}_n^{k}$ be an arbitrary unicyclic graph. Then $F(G_{k ,2}^{(n)}) > F(G_{k + 1 ,2}^{(n)})$.
\end{thm}

\n\textit{Proof.}
From the definition of F-index, we have
\begin{eqnarray*}
F(G_{k ,2}^{(n)}) - F(G_{k + 1 ,2}^{(n)}) &=& (n - k + 1)^3 + 1^3 - 2^3 - (n - k)^3\\
&=& 3(n - k)^2 + 3(n - k) - 7 > 0 ~\mbox{since} ~n - k \geq 2.
\end{eqnarray*}

Hence, $F(G_{k ,2}^{(n)}) > F(G_{k + 1, 2}^{(n)})$.
\qed

\begin{thm}\label{43}
Let $p \geq q \geq r \geq 0$, and $p + q +r = n - 3$. Then \\
(i) $F(S_{p, q, r}) < F(S_{p + 1, q - 1, r})$, and
(ii) $F(S_{p, q, r}) < F(S_{p, q + 1, r - 1})$.
\end{thm}

\n\textit{Proof.} (i) By the definition of F-index, we have
\begin{eqnarray*}
F(S_{p, q, r}) - F(S_{p + 1, q - 1, r}) &=&  (p + 2)^3 (q + 2)^3 (r + 2)^3 - (p + 3)^3 - (q + 1)^3 - (r + 2)^3\\
&=& (q - p -1)(3p + 3q + 12) < 0.
\end{eqnarray*}
Hence, $ F(S_{p, q, r}) < F(S_{p + 1, q - 1, r})$.

\noindent{(ii)} \begin{eqnarray*}
F(S_{p, q, r}) - F(S_{p, q + 1, r - 1}) &=& (p + 2)^3 + (q + 2)^3 + (r + 2)^3 - (p + 2)^3 - (q + 3)^3 - (r + 1)^3\\
&=& (r - q - 1)(3r + 3q + 12) < 0.
\end{eqnarray*}
Hence, $ F(S_{p, q, r}) < F(S_{p, q + 1, r - 1})$.
\qed

By attaching $K_{1, n - 5}$ and $K_{1, 2}$ to the adjacent vertices of $C_{3}$, respectively, we have the graph $G_{3 ,3}^{(n)}$. Also by attaching $K_{1, n - 5}$ to one vertex of $C_{3}$ and $K_{1, 2}$ to another two vertices of $C_{3}$, respectively, we have the graph $G_{3 ,4}^{(n)}$. See Figure 8(a) and 8(b) for illustration.

 Now, $F(G_{3, 3}^{(n)}) = (n - 3)^3 + 2^3 + 4^3 + (n - 2)1^3 = n^3 - 9n^2 + 28n + 43$,\\
 and $F(G_{3, 4}^{(n)}) = (n - 3)^3 + 3^3 + 3^3 + (n - 2)1^3 = n^3 - 9n^2 + 28n + 25$.
 
By Lemma $5$, Lemma $6$, Theorem $3$ and above calculation, we have the following theorem.

\begin{thm}\label{44}
For $ n \geq 6$, we have $F(G_{3, 3}^{(n)}) > F(G_{3, 4}^{(n)}) > . . .$
\end{thm}

Let $\mathcal{S}$ denote the set of graphs $S_{p, q, r}$. Then by Lemma $5$, Lemma $6$, Theorem $3$, and  Theorem $4$ we have the following theorem.

\begin{thm}\label{45}
For $n \geq 6$, the order in $\mathcal{S}$ with respect to the F-indices is $F(G_{3, 1}^{(n)}) > F(G_{3, 2}^{(n)}) > F(G_{3, 3}^{(n)}) > F(G_{3, 4}^{(n)}) > . . .$
\end{thm}

Let $\mathcal{U}_n^{'3}$ be the set of unicyclic graphs with $C_3$ being the only cycle and there is at least one vertex in the graph, which is at distant $\geq 2$  from $C_{3}$. Obviously, $\mathcal{U}_n^{'3} = \mathcal{U}_n^{3} - \mathcal{S}$. By Lemma 1 and Lemma 2, it is clear that the graphs with the largest F-indices in $\mathcal{U}_n^{'3}$ must be made by attaching $K_{1, l}(l \geq 1)$ to one of the pendent vertices of $S_{i, j, k}(i, j \geq 0, k \geq 1)$. Denote the graph as $R_{i, j, k, l}$, as shown in Figure 9.

Similar to Theorem $3$, we have the following theorem.

\begin{thm}\label{46}
Let $i \geq j \geq 1, i + j + k + l = n - 3$. Then, we have $F(R_{i + 1, j - 1, k, l}) > F(R_{i, j, k, l})$. In particular, $F(R_{i + j, 0, k, l}) > F(R_{i, j, k, l}).$
\end{thm}

\n\textit{Proof.} From the definition of F-index, we have
\begin{eqnarray*}
F(R_{i + 1, j - 1, k, l}) &=& (i + 3)^3 + (j + 1)^3 + (k + 2)^3 + (l + 1)^3 + (n - 4),\\
\mbox{and } ~F(R_{i, j, k, l}) &=& (i + 2)^3 + (j + 2)^3 + (k + 2)^3 + (l + 1)^3 + n - 4.\\
 \mbox{So, } ~F(R_{i + 1, j - 1, k, l}) - F(R_{i, j, k, l}) &=& (i + 3)^3 + (j + 1)^3 - (i + 2)^3 - (j + 2)^3\\
&=& 3(i^2 - j^2 + 5i - 3j + 4) > 0.
\end{eqnarray*}

Hence, $F(R_{i + 1, j - 1, k, l}) > F(R_{i, j, k, l})$.

In particular,

$F(R_{i + j, 0, k, l}) > F(R_{i, j, k, l})$.
\qed
\begin{thm}\label{47}
Let $i \geq 1$. Then we have $F(R_{i, 0, k, l}) < F(R_{0, 0, i + k, l})$.
\end{thm}

\n\textit{Proof.}
From definition of F-index, we have
\begin{eqnarray*}
F(R_{i, 0, k, l}) &=& (i + 2)^3 + (k + 2)^3 + (l + 2)^3 + (l + k + i - 1)1^3,\\
 \mbox{and } ~F(R_{0, 0, i + k, l}) &=& 2^3 + 2^3 + (i + k + 2)^3 + (l + 2)^3 + (l + k + i - 1).\\
 \mbox{Then } ~F(R_{0, 0, i + k, l}) - F(R_{0, 0, k, l}) &=& 16 +  (i + k + 2)^3 - (i + 2)^3 - (k + 2)^3 - 2^3\\
&=& 3ik(i + k + 4) > 0.
\end{eqnarray*}
Hence, $F(R_{i, 0, k, l}) < F(R_{0, 0, i + k, l})$.
\qed

We shall denote $R_{0, 0, k, 1}$ simply as $R_{k, 1}$, and which is also denoted by $G_{3, 1}^{'(n)}$.

\begin{thm}\label{48}
For $ n \geq 9$, we have $F(G_{3, 3}^{(n)}) < F(R_{k, 1}) < F(G_{3, 2}^{(n)})$.
\end{thm}

\n\textit{Proof.}
By definition of F-index, we have
\begin{eqnarray*}
F(G_{3, 3}^{(n)}) &=& n^3 - 9n^2 + 28n + 43,\\
F(G_{3, 2}^{(n)}) &=& n^3 - 6n^2 + 13n + 24,\\
\mbox{and} ~F(R_{k, 1}) &=& n^3 - 6n^2 + 13n + 12\\
\mbox{Therefore}, ~F(R_{k, 1}) - F(G_{3, 3}^{(n)})&=& 3n^2 - 15n - 13 > 0,\\
\mbox{and} ~F(G_{3, 2}^{(n)}) - F(R_{k, 1}) &=& 12 > 0.
\end{eqnarray*}
Hence, $F(G_{3, 3}^{(n)}) < F(R_{k, 1}) < F(G_{3, 2}^{(n)})$.
\qed

\begin{thm}\label{49}
For $ n \geq 6$, $l \geq 2$, $k + l + 3 = n$, we have $F(G_{3, 4}^{(n)}) > F(R_{k, l})$.
\end{thm}

\n\textit{Proof.}
From the definition of $R_{k, l}$ and since $ k + l + 3 = n$, we have
\begin{eqnarray*}
F(R_{k, l}) &=& 2^3 + 2^3 + (k + 2)^3 + 2^3 + k\\
&=& 24 + (n - l - 1)^3 + (n - 3 - l)\\
&=& n^3 - 3n^2 + 4n + 12 + 3nl^2 + (6n - 3n^2)l.
\end{eqnarray*}
Let $f(l) = n^3 - 3n^2 + 4n + 12 + 3nl^2 + (6n - 3n^2)l$, $l \in [2, n - 4]$.\\
Then $\max {f(l)} = \{f(2), f(n - 4)\} = n^3 - 9n^2 + 28n + 12$, (since $f(2) = f(n - 4)$).\\
Therefore, $F(G_{3, 4}^{(n)}) > F(R_{k, l})$.
\qed

\begin{thm}\label{410}
Let $ n \geq 9$. Then F-indices order in $\mathcal{U}_n^{3}$ is\\
$F(G_{3, 1}^{(n)}) > F(G_{3, 2}^{(n)}) > F(G_{3, 1}^{'(n)}) > F(G_{3, 3}^{(n)}) > F(G_{3, 4}^{(n)}) > F(R_{k, l}) > ...$
\end{thm}

Let $G_{4, 3}^{(n)}$ be the graph obtained from a $C_{4}$ by attaching $n - 5$ leaves to one of its vertices and another leaf to the vertex which is at a distance 2 from the $n - 3$ degree vertex of $C_{4}$.

By the definition of F-index, we have

$F(G_{4, 3}^{(n)}) = n^3 - 9n^2 + 28n + 12$,\\ and
$F(G_{4, 2}^{(n)}) = n^3 - 9n^2 + 28n + 12$.

Similar to Theorem 10, we have the following theorem.

\begin{thm}\label{411}
Let $ n \geq 6$, then F-index order in $\mathcal{U}_n^{4}$ is\\
$F(G_{4, 1}^{(n)}) > F(G_{4, 2}^{(n)}) = F(G_{4, 3}^{(n)}) > F(G_{4, 1}^{'(n)}) > ...$, where $G_{4, 1}^{'(n)}$ is obtained from $G_{4, 1}^{(n - 1)}$ by attaching $K_{2}$ to one pendent edges of $G_{4, 1}^{(n - 1)}$.
\end{thm}

\begin{thm}\label{412}
Let $n \geq 6$. Then we have\\
(i) $F(G_{4, 1}^{(n)}) = F(G_{3, 1}^{'(n)})$, and\\
(ii) $F(G_{3, 4}^{(n)}) > F(G_{4, 2}^{(n)}).$
\end{thm}

\n\textit{Proof.} (i) By the definition of F-index, we have\\
$F(G_{4, 1}^{(n)}) = 2^3 + 2^3 + 2^3 + (n - 2)^3 + (n - 4) = n^3 - 6n^2 + 13n + 12 = F(G_{3, 1}^{'(n)})$.\\
\noindent (ii) $F(G_{3, 4}^{(n)}) - F(G_{4, 2}^{(n)}) = (n^3 - 9n^2 + 28n + 25) - (n^3 - 9n^2 + 28n + 12) = 13$.\\
Hence, $F(G_{3, 4}^{(n)}) > F(G_{4, 2}^{(n)})$.
\qed

Combining all the above results, we have the upper bounds of unicyclic graphs with respect to F-indices.

\begin{thm}\label{413}
For $n \geq 6$, we have\\
$F(G_{3, 1}^{(n)}) > F(G_{3, 2}^{(n)}) > F(G_{3, 1}^{'(n)}) = F(G_{4, 1}^{(n)}) > F(G_{3, 3}^{(n)}) > F(G_{3, 4}^{(n)}) > ...$
\end{thm}

\section{The lower bounds of the unicyclic graphs with respect to F-indices}

Given integers $n$ and $k$ with $3 \leq k \leq n - 1$, the lollipop $L_{n, k}$ is the unicyclic graph of order $n$ obtained from the two vertex disjoint graphs $C_{k}$ and $P_{n -k}$ by adding an edge joining a vertex of $C_{k}$ to an end vertex of $P_{n - k}$.

\begin{thm}\label{51}
The graph $C_{n}$ is the unique graph with the smallest F-index among all unicyclic graphs with $n$ vertices.
\end{thm}

\n\textit{Proof.}
First we shall prove that if $G$ is a unicyclic graph with $n \geq 7$ vertices, then $F(G)$ attains the smallest value only if degree sequence of $G$ is $[2^n]$. Suppose $F(G)$ attains the smallest value and degree sequence of $G$ is not equal to $[2^n]$. Let $C = u_{1}u_{2}. . . u_{k}u_{1}$ be the unique cycle in $G$. Then $k < n$ and there is at least one vertex $u_{i}$ with $d(u_{i}) \geq 3$. Without loss of generality, we assume $d(u_{1}) \geq 3$. Choose a maximal $C-$path $P[u_{1}, v_{1}]$ in $G$. Clearly $d(v_{1}) = 1$. Let $G^{'} = G - u_{1}u_{2} + u_{2}v_{1}$. Then we have, $F(G) - F(G^{'}) = [d^3(u_{1}) + d^3(v_{1})] - [(d^3(u_{1}) - 1)^3 + (d^3(v_{1}) + 1)^3] = 3(d(u_{1}) + d(v_{1}))(d(u_{1}) - d(v_{1}) - 1) > 0$. Therefore $F(G) > F(G^{'})$, a contradiction. Hence degree sequence of $G$ is $[2^n]$.

Since among all unicyclic graph, the cycle $C_{n}$ has only degree sequence $[2^n]$, it has smallest F-index among all unicyclic graph with $n$ vertices.
\qed

\begin{thm}\label{52}
Let $G \in \mathcal{U}_n^{k}$, $3 \leq k \leq n - 1$ be an arbitrary unicyclic graph. Then $F(G) \geq F(L_{n, k})$, with equality if and only if $G \cong L_{n, k}$.
\end{thm}

\n\textit{Proof.}
By Transformation C, D and Lemma 3 and 4, the conclusion is obvious.
\qed

\begin{thm}\label{53}
Let $G \in \mathcal{U}_n - C_{n}$ be an arbitrary unicyclic graph, then $F(G) > F(L_{n, k})$, $k \in \{3, 4, ..., n - 1\}$.
\end{thm}

\n\textit{Proof.}
By the definition of $L_{n, k}$, we have $F(L_{n, k}) = 8n + 14$. Therefore value of $F$ are function of $n$, not related to $k$, hence $F(L_{n, k}) = F(L_{n, l})$, where $k \in \{3, 4, . . ., n - 1\}$ and $k \neq l$, i.e., if $G$ is not isomorphic to $C_{n}$, then $F(G) > F(L_{n, k})$ for $k \in \{3, 4, ..., n - 1\}.$
\qed

\section*{Acknowledgement}
This work has been partially supported by University Grants Commission, India, through Grant no. F1-17.1/2016-17/MANF-2015-17-WES-60163, awarded to the first author, under Maulana Azad National Fellowship Scheme.


\begin{figure}[h]
\begin{center}
\includegraphics[width=0.75\textwidth]{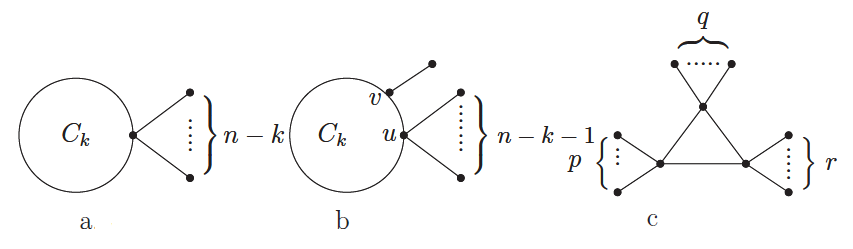}
\caption{(a) $G_{k, 1}^{(n)}$  (b) $G_{k, 2}^{(n)}$ (c) $S_{p, q, r}$. }
\end{center}
\end{figure}
\begin{figure}[h]
\begin{center}
\includegraphics[width=0.75\textwidth]{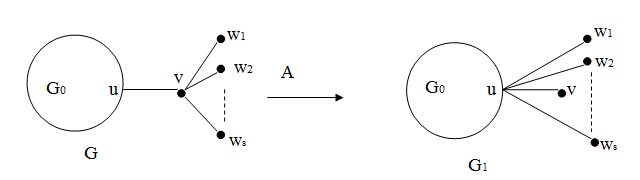}
\caption{Transformation A }
\end{center}
\end{figure}
\begin{figure}[h]
\begin{center}
\includegraphics[width=0.65\textwidth]{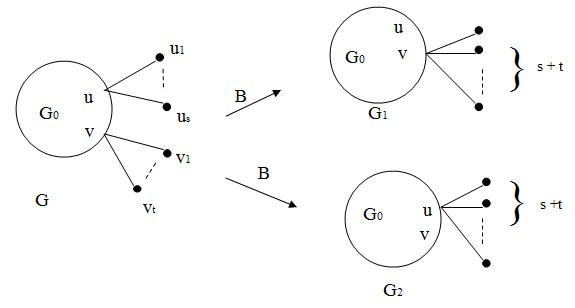}
\caption{Transformation B }
\end{center}
\end{figure}

\begin{figure}[h]
\begin{center}
\includegraphics[width=0.75\textwidth]{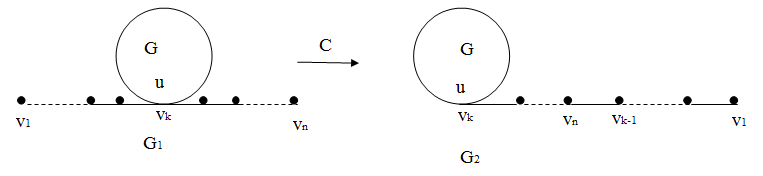}
\caption{Transformation C }
\end{center}
\end{figure}
\begin{figure}[h]
\begin{center}
\includegraphics[width=0.75\textwidth]{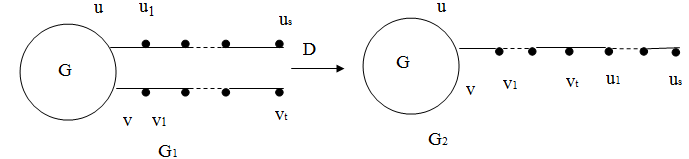}
\caption{Transformation D }
\end{center}
\end{figure}
\begin{figure}[h]
\begin{center}
\includegraphics[width=0.75\textwidth]{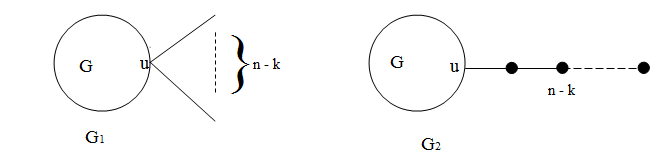}
\caption{Transformation E }
\end{center}
\end{figure}
\begin{figure}[h]
\begin{center}
\includegraphics[width=0.75\textwidth]{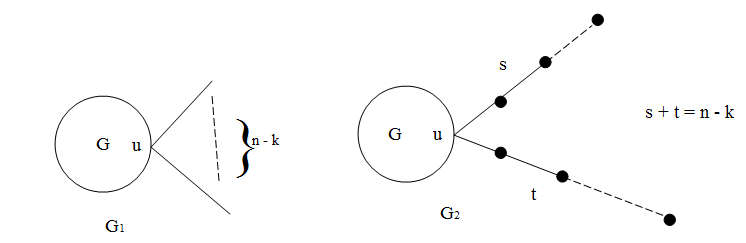}
\caption{Transformation F }
\end{center}
\end{figure}
\begin{figure}[h]
\begin{center}
\includegraphics[width=0.75\textwidth]{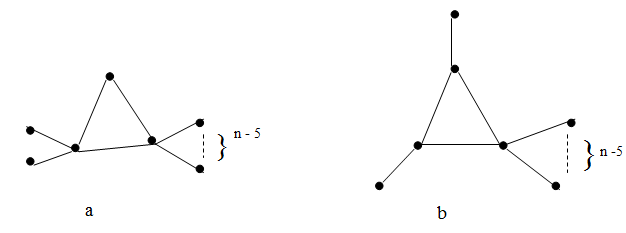}
\caption{(a) $G_{3, 3}^{(n)}$. (b) $G_{3, 3}^{(n)}$. }
\end{center}
\end{figure}
\begin{figure}[h]
\begin{center}
\includegraphics[width=0.75\textwidth]{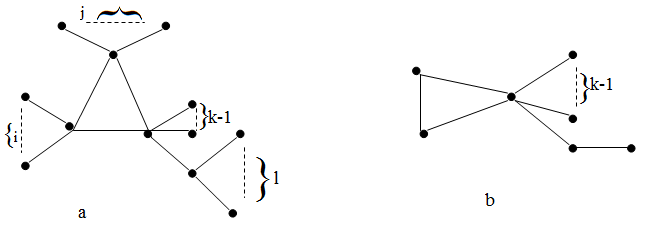}
\caption{(a) $R_{i, j, k, l}$. (b) $R_{0, 0, k, 1}$. }
\end{center}
\end{figure}

\end{document}